# FUNCTIONAL EQUATIONS OF TORNHEIM ZETA FUNCTION

## AND THE α-CALCULUS OF $\dfrac{\partial^r}{\partial s^r}\zeta(s,\alpha)$.

V.V.RANE

A-3/203, ANAND NAGAR ,

DAHISAR ,MUMBAI-400 068,

INDIA

v_v_rane@yahoo.co.in

**Abstract :** For the Tornheim double zeta function $T(s_1,s_2,s_3)$ of complex variables, we obtain its functional equations , which are new . Using the calculus of $\dfrac{\partial^r}{\partial s^r}\zeta(s,\alpha)$ as a function of α( developed in author[7] ) as the tool, we obtain with ease an expression for $T(s,s,s)$ for any complex s ; and evaluate $T(n,n,n)$ for any integer $n \geq 1$ ; and also evaluate $T(n_1,0,n_2)$ for integers $n_1,n_2 > 1$ with $n_1+n_2$ odd .

We obtain for Tornheim zeta function , the counterparts of Euler's formula for $\zeta(2n)$ and its analogue for $\zeta(2n+1)$, where $\zeta(s)$ is Riemann zeta function .

# FUNCTIONAL EQUATIONS OF TORNHEIM ZETA FUNCTION

AND THE α-CALCULUS OF $\dfrac{\partial^r}{\partial s^r} \zeta(s,\alpha)$.

V.V.RANE

A-3/203, ANAND NAGAR,

DAHISAR, MUMBAI-400 068,

INDIA

v_v_rane@yahoo.co.in

For complex $\alpha \neq 0, -1, -2, \ldots$ and for the complex variable s, let

$\zeta(s,\alpha)$ be the Hurwitz zeta function defined by $\zeta(s,\alpha) = \sum\limits_{n \geq 0} (n+\alpha)^{-s}$ for Re s>1 ; and its

analytic continuation. Let $\zeta(s,1) = \zeta(s)$, the Riemann zeta function. For an integer $r \geq 0$, let $\zeta^{(r)}(s,\alpha) = \dfrac{\partial^r}{\partial s^r} \zeta(s,\alpha)$. For the complex variables $s_1, s_2, s_3$, let Tornheim's

double zeta function $T(s_1, s_2, s_3)$ be defined by

$T(s_1,s_2,s_3) = \sum\limits_{n_1 \geq 1} \sum\limits_{n_2 \geq 1} n_1^{-s_1} n_2^{-s_2} \cdot (n_1+n_2)^{-s_3}$ for Re$(s_2+s_3)$, Re$(s_3+s_1) > 1$ and

Re $(s_1+s_2+s_3) > 2$; and its analytic continuation. Then it is known that the

function $T(s_1,s_2,s_3)$ can be continued as a meromorphic function to the

3-dimensional complex field $\mathbb{C}^3$.

For an integer $r \geq 1$, we define the multiple zeta function by

$$\zeta_r(s_1, s_2, \ldots, s_r) = \sum_{n_1 \geq 1} n_1^{-s_1} \cdot \sum_{n_2 > n_1} n_2^{-s_2} \cdot \sum_{n_3 > n_2} n_3^{-s_3} \ldots \cdot \sum_{n_r > n_{r-1}} n_r^{-s_r}$$

for $\mathrm{Re}(s_1 + s_2 + \ldots + s_r) > r$ and its analytic continuation.

We make following observations

1) $T(s_1, s_2, s_3) = T(s_2, s_1, s_3)$

2) $T(s_1, s_2, 0) = \zeta(s_1) \cdot \zeta(s_2)$

3) $T(s_1, 0, s_2) = \zeta_2(s_1, s_2)$

4) $\zeta_2(s_1, s_2) + \zeta_2(s_2, s_1) = \zeta(s_1) \cdot \zeta(s_2) - \zeta(s_1 + s_2)$

so that $T(s_1, 0, s_2) + T(s_2, 0, s_1) = \zeta(s_1) \cdot \zeta(s_2) - \zeta(s_1 + s_2)$

5) $T(s_1, s_2 - 1, s_3 + 1) + T(s_1 - 1, s_2, s_3 + 1) = T(s_1, s_2, s_3)$.

For an integer $r \geq 0$, we shall write $\zeta^{(r)}(s, \alpha) = \dfrac{\partial^r}{\partial s^r} \zeta(s, \alpha)$. In what follows $\Gamma(s)$ shall stand for Euler's gamma function. We define Bernoulli polynomial $B_n(\alpha)$ (of degree n) in variable α, by $\dfrac{z e^{\alpha z}}{e^z - 1} = \sum_{n \geq 0} \dfrac{B_n(\alpha)}{n!} z^n$ for $|z| < 2\pi$.

Note that $B_0(\alpha) = 1$. We write $B_n = B_n(0)$. Then $B_n$'s are called Bernoulli numbers, which are rational numbers. We also note that $\zeta(-n, \alpha) = -\dfrac{B_{n+1}(\alpha)}{n+1}$ for integral $n \geq 0$. We also have $B_n(1 - \alpha) = (-1)^n B_n(\alpha)$.

Using the theory of $\zeta^{(r)}(s, \alpha)$ as an analytic function of the complex variable α in author [5],[6], we have developed the calculus and analysis of $\zeta^{(r)}(s, \alpha)$ as a function of α in author [7]. Using the calculus of $\zeta^{(r)}(s, \alpha)$ as a function of α  (of author [7]) , and the calculus as developed in the present paper , we obtain  functional equations of the Tornheim double zeta function $T(s_1,s_2,s_3)$ , which are  the counterparts of the functional equation of Riemann's zeta function, which are new . Using our version of Hurwitz's  formula for  $\zeta(s,\alpha)$ , we have developed  functional equations of $T(s_1,s_2,s_3)$ . Actually , our method is capable of giving  functional equations of the r-fold Tornheim zeta function .For complex $\alpha \neq 0,-1,-2,\ldots\ldots\ldots\ldots$ and for real λ  with 0<λ≤1 , define  Lerch's zeta function $\phi(\lambda,\alpha,s)$ by $\phi(\lambda,\alpha,s) = \sum_{n\geq 0} e^{2\pi i \lambda n}(n+\alpha)^{-s}$ for Re s>1 ; and its analytic continuation .If we use in the place( of our version) of Hurwitz's  formula for $\zeta(s,\alpha)$ , the fact that for Re s<1 , 0<λ<1and 0<α<1 ,

$$\phi(\lambda,\alpha,s) = i(2\pi)^{s-1}\Gamma(1-s)e^{-\frac{\pi i s}{2}}\sum_{|n|\geq 0} e^{-2\pi i(n+\lambda)\alpha}\cdot(n+\lambda)^{s-1} = \Gamma(1-s)\sum_{|n|\geq 0} e^{2\pi i(n-\lambda)\alpha}\left(2\pi i(n-\lambda)\right)^{s-1},$$

we get  functional equations for Tornheim –Lerch's double zeta function .

Our results have been obtained in terms of integrals of the type $\int_0^1 \zeta(s_1,\alpha)\zeta(s_2,\alpha)\zeta^{(r)}(s_3,\alpha)d\alpha$ , where $s_1,s_2,s_3$, are complex variables and $r\geq 0$ is an integer .Espinosa and Moll [1] have also provided an explicit formula for T($s_1,s_2,s_3$) for real $s_1,s_2,s_3$  in terms of the integrals

$$I_1 = \int_0^1 \zeta(1-s_1,1-\alpha)\zeta(1-s_2,\alpha)\zeta(1-s_3,\alpha)d\alpha, \ I_2 = \int_0^1 \zeta(1-s_1,\alpha)\zeta(1-s_2,1-\alpha)\zeta(1-s_3,\alpha)d\alpha,$$

$$I_3 = \int_0^1 \zeta(1-s_1,\alpha)\zeta(1-s_2,\alpha)\zeta(1-s_3,1-\alpha)d\alpha \text{ and } I_4 = \int_0^1 \zeta(1-s_1,\alpha)\zeta(1-s_2,\alpha)\zeta(1-s_3,\alpha)d\alpha.$$

Note that as the functions of complex variables $s_1,s_2,s_3$, $I_1,I_2,I_3,I_4$ are analytic functions of $s_1, s_2, s_3$ for Re $s_1>1$, Re $s_2>1$, Re $s_3>1$ respectively, in view of the remark following Lemma 2 below.

We note that for complex numbers $s_1, s_2, s_3$,

1) $\int_0^1 \zeta(1-s_1,\alpha)\zeta(1-s_2,\alpha)\zeta(1-s_3,\alpha)d\alpha = \int_0^1 \zeta(1-s_1,1-\alpha)\zeta(1-s_2,1-\alpha)\cdot\zeta(1-s_3,1-\alpha)d\alpha$

2) $\int_0^1 \zeta(1-s_1,1-\alpha)\zeta(1-s_2,1-\alpha)\cdot\zeta(1-s_3,\alpha)d\alpha = \int_0^1 \zeta(1-s_1,\alpha)\zeta(1-s_2,\alpha)\cdot\zeta(1-s_3,1-\alpha)d\alpha$

3) If $m_1 \geq 1$ is an integer and $s_2, s_3$ complex, then

$$\int_0^1 \zeta(1-m_1,1-\alpha)\zeta(1-s_2,\alpha)\zeta(1-s_3,\alpha)d\alpha = -\frac{1}{m_1}\int_0^1 B_{m_1}(1-\alpha)\zeta(1-s_2,\alpha)\cdot\zeta(1-s_3,\alpha)d\alpha$$

$$= -\frac{(-1)^{m_1}}{m_1}\int_0^1 B_{m_1}(\alpha)\zeta(1-s_2,\alpha)\cdot\zeta(1-s_3,\alpha)d\alpha = (-1)^{m_1}\int_0^1 \zeta(1-m_1,\alpha)\zeta(1-s_2,\alpha)\cdot\zeta(1-s_3,\alpha)d\alpha$$

4) For integers $m_1, m_2, m_3 \geq 1$, $\int_0^1 \zeta(1-m_1,\alpha)\cdot\zeta(1-m_2,\alpha)\cdot\zeta(1-m_3,\alpha)d\alpha$ is a rational number, which equals zero, if N= $m_1 + m_2 + m_3$ is odd. (See Lemma 6 below.)

5) If $m_1, m_2 \geq 1$ are integers and $s$ is complex and $r \geq 0$ an integer, then

$$\int_0^1 \zeta(1-m_1,\alpha)\zeta(1-m_2,\alpha)\zeta^{(r)}(1-s,1-\alpha)d\alpha = \int_0^1 \zeta(1-m_1,1-\alpha)\zeta(1-m_2,1-\alpha)\cdot\zeta^{(r)}(1-s,\alpha)d\alpha$$

$$= \frac{1}{m_1 m_2}\int_0^1 B_{m_1}(1-\alpha)B_{m_2}(1-\alpha)\zeta^{(r)}(1-s,\alpha)d\alpha = \frac{(-1)^{m_1+m_2}}{m_1 m_2}\int_0^1 B_{m_1}(\alpha)B_{m_2}(\alpha)\zeta^{(r)}(1-s,\alpha)d\alpha$$

$$= (-1)^{m_1+m_2}\int_0^1 \zeta(1-m_1,\alpha)\zeta(1-m_2,\alpha)\cdot\zeta^{(r)}(1-s,\alpha)d\alpha.$$

We have indicated in author [7] that if $m_1, m_2 \geq 1$ are integers and Re s<1 , then $\int_0^1 \zeta(1-m_1,\alpha)\cdot\zeta(1-m_2,\alpha)\cdot\zeta^{(r)}(s,\alpha)d\alpha$ is explicitly computable as a linear combination of $\zeta^{(\ell)}(s-1), \zeta^{(\ell)}(s-2)............\zeta^{(\ell)}(s-N))$ for $0\leq \ell \leq r$ with coefficients dependent on s, where N=degree of the polynomial

$$\zeta(1-m_1,\alpha)\cdot\zeta(1-m_2,\alpha) = \frac{1}{m_1 m_2} B_{m_1}(\alpha)\cdot B_{m_2}(\alpha)\ .$$ Thus we see that

$\int_0^1 \zeta(1-n_1,\alpha)\cdot\zeta(1-n_2,\alpha)\cdot\zeta'(1-n_3,\alpha)d\alpha$ is a linear combination of

$\zeta'(-n_3), \zeta'(-n_3-1), \zeta'(-n_3-2),............,\zeta'(-n_1-n_2-n_3)$ with rational coefficients .

Consequently , we find that $(-1)^{n_3}T(n_1,n_2,n_3)+(-1)^{n_1}T(n_2,n_3,n_1)+(-1)^{n_2}T(n_3,n_1,n_2)$

is a rational multiple of $\pi^N$ for $N=n_1+n_2+n_3$ even ; and

$(-1)^{n_3}T(n_1,n_2,n_3)+(-1)^{n_1-1}\cdot T(n_2,n_3,n_1)+(-1)^{n_2-1}\cdot T(n_3,n_1,n_2)$ equals $\pi^{N-1}$ multiplied by

a linear combination of $\zeta'(-n_3), \zeta'(-n_3-1), \zeta'(-n_3-2),............,\zeta'(-n_1-n_2-n_3)$ with rational coefficients , for $N=n_1+n_2+n_3$ odd . These two facts about Tornheim double zeta function are the counterparts of the two facts about Riemann zeta function namely , $\zeta(N) = \dfrac{(-1)^{\frac{N}{2}} \pi^N 2^{N-1}}{(N-1)!}\cdot\zeta(1-N)$ for $N\geq 2$ even ;

and $\zeta(N) = \dfrac{(-1)^{\frac{N-1}{2}}\cdot 2^N \pi^{N-1}}{(N-1)!}\cdot\zeta'(1-N)$ for $N\geq 3$ odd . (See author [5]).

Note that Euler's formula for $\zeta(N)$ for $N\geq 2$ even , is actually the functional equation for $\zeta(s)$ at positive even integer argument .

The above equation for $\zeta(N)$ for $N\geq 3$ odd ( of author [5]) is actually the counterpart at odd positive integer argument , of Euler's formula for $\zeta(s)$ (at even positive integer argument ).

Next we state our Theorem .

**Theorem :** Let $s_1, s_2, s_3$ be complex numbers and let $N = s_1 + s_2 + s_3$.

Then we have

I) $\int_0^1 \zeta(1-s_1,\alpha)\zeta(1-s_2,\alpha)\zeta(1-s_3,\alpha)d\alpha = 2(2\pi)^{-(s_1+s_2+s_3)} \cdot \Gamma(s_1)\Gamma(s_2)\Gamma(s_3) \cdot$

$\cdot \left\{ \cos\frac{\pi}{2}(N-2s_3) \cdot T(s_1,s_2,s_3) + \cos\frac{\pi}{2}(N-2s_1) \cdot T(s_2,s_3,s_1) + \cos\frac{\pi}{2}(N-2s_2) \cdot T(s_3,s_1,s_2) \right\}$

II) $\int_0^1 \zeta(1-s_1,\alpha)\zeta(1-s_2,\alpha)\zeta(1-s_3,1-\alpha)d\alpha = 2(2\pi)^{-(s_1+s_2+s_3)} \cdot \Gamma(s_1)\Gamma(s_2)\Gamma(s_3) \cdot$

$\cdot \left\{ \cos\frac{\pi}{2}N \cdot T(s_1,s_2,s_3) + \cos\frac{\pi}{2}(N-2s_2) \cdot T(s_2,s_3,s_1) + \cos\frac{\pi}{2}(N-2s_1) \cdot T(s_3,s_1,s_2) \right\}$

Remark : 1) Note that in general , there are two more functional equations

corresponding to the integrals

$I_1 = \int_0^1 \zeta(1-s_1,1-\alpha)\zeta(1-s_2,\alpha) \cdot \zeta(1-s_3,\alpha)d\alpha, \quad I_2 = \int_0^1 \zeta(1-s_1,\alpha)\zeta(1-s_2,1-\alpha)\zeta(1-s_3,\alpha)d\alpha.$

**Corollary 1:** Let $n_1, n_2 \geq 1$ be integers and s complex . Let $N = n_1 + n_2 + s$.

Then we have

a) $\int_0^1 \zeta(1-n_1,\alpha) \cdot \zeta(1-n_2,\alpha) \cdot \zeta(1-s,\alpha)d\alpha = 2(2\pi)^{-N} \cdot \Gamma(n_1)\Gamma(n_2) \cdot \Gamma(s) \cdot$

$\cdot \left\{ \cos\frac{\pi}{2}(N-2s)T(n_1,n_2,s) + (-1)^{n_1} \cos\frac{\pi N}{2} \cdot T(n_2,s,n_1) + (-1)^{n_2} \cos\frac{\pi N}{2} \cdot T(s,n_1,n_2) \right\}$

b) $\int_0^1 \zeta(1-n_1,\alpha)\cdot\zeta(1-n_2,\alpha)\cdot\zeta(1-s,1-\alpha)d\alpha = 2(2\pi)^{-N}\Gamma(n_1)\Gamma(n_2)\Gamma(s)\cdot\cos\frac{\pi}{2}N \cdot$

$\cdot\{T(n_1,n_2,s)+(-1)^{n_2}\cdot T(n_2,s,n_1)+(-1)^{n_1}\cdot T(s,n_1,n_2)\}$

**Remark :** The above two equations of Corollary 1 are actually identical **in view of**

**the fact that** $\int_0^1 \zeta(1-n_1,\alpha)\zeta(1-n_2,\alpha)\zeta(1-s,1-\alpha)d\alpha$

$= (-1)^{n_1+n_2}\int_0^1 \zeta(1-n_1,\alpha)\zeta(1-n_2,\alpha)\zeta(1-s,\alpha)d\alpha.$

Corollary 1 b) is our version for the functional relation given by

Tsumura [8] in the form

$T(n_1,n_2,s)+(-1)^{n_2}T(n_2,s,n_1)+(-1)^{n_1}T(s,n_1,n_2) = N_0(n_1,n_2,s)$, where $N_0(n_1,n_2,s)$ has a

complicated explicit expression . Nakamura [2] has given a simpler alternative expression for $N_0(n_1,n_2,s)$ and a simpler proof of Tsumura's functional relation .

Here it is to be emphasized that $\int_0^1 \zeta(1-n_1,\alpha)\cdot\zeta(1-n_2,\alpha)\zeta(1-s,1-\alpha)d\alpha$ is explicitly

computable as a linear combination of $\zeta(-s),\zeta(-s-1),\ldots,\zeta(-s-M)$ (with

coefficients dependent upon s ), where M= degree of the polynomial

$\zeta(1-n_1,\alpha)\cdot\zeta(1-n_2,\alpha)$ .

**Corollary 2 :** Let $n_1,n_2,n_3 \geq 1$ be integers and let $N = n_1+n_2+n_3$ . Then we have

$\int_0^1 \zeta(1-n_1,\alpha)\cdot\zeta(1-n_2,\alpha)\zeta(1-n_3,\alpha)d\alpha = 2(2\pi)^{-N}\cdot\Gamma(n_1)\Gamma(n_2)\Gamma(n_3)\cos\frac{\pi N}{2} \cdot$

$\cdot\{(-1)^{n_3}T(n_1,n_2,n_3)+(-1)^{n_1}T(n_2,n_3,n_1)+(-1)^{n_2}T(n_3,n_1,n_2)\}$

**Note :** If $N = n_1 + n_2 + n_3$ is odd, then either side of the equation equals zero .

From the statement of our Theorem , it is clear that we have four equations in three unknowns $T_1 = T(s_1, s_2, s_3)$, $T_2 = T(s_2, s_3, s_1)$ and $T_3 = T(s_3, s_1, s_2)$.

Writing $d_i = d_i(s_1, s_2, s_3) = \dfrac{(2\pi)^{s_1+s_2+s_3}}{2\Gamma(s_1)\Gamma(s_2)\Gamma(s_3)} I_i$ for $i = 1,2,3,4$,

we have four equations of the type

$a_1T_1+b_1T_2+c_1T_3=d_1$ ; $a_2T_1+b_2T_2+c_2T_3=d_2$ ; $a_3T_1+b_3T_2+c_3T_3=d_3$ ; $a_4T_1+b_4T_2+c_4T_3=d_4$,

where $a_i=a_i(s_1,s_2,s_3)$ , $b_i=b_i(s_1,s_2,s_3)$ , $c_i=c_i(s_1,s_2,s_3)$ are known smooth functions of $s_1,s_2,s_3$ for i=1,2,3 ,4 and $d_i=d_i(s_1,s_2,s_3)$ for i=1,2,3,4 are smooth functions of $s_1,s_2,s_3$ . Solving any three equations in three unknowns $T_1,T_2,T_3$, we get the values of $T_1,T_2,T_3$, for general complex numbers $s_1,s_2,s_3$ . When $s_1=n_1\geq 1$ is an integer and $s_2,s_3$ are complex , we are left with only two equations in view of the observations made above in respect of integrals with respect to α . However when $n_1,n_2 \geq 1$ are integers and s complex , all the four equations are identical , giving a single equation $a'T(n_1,n_2,s)+b'T(n_2,s,n_1)+c'T(s,n_1,n_2)=e'$, where a',b',c',e' are known smooth functions of $n_1,n_2,s$ . In particular , when $n_1,n_2,n_3$ are non-negative integers , all the four equations coincide to give a single equation of the form a $T(n_1,n_2,n_3)$+b $T(n_2,n_3,n_1)$+c $T(n_3,n_1,n_2)=d\pi^N$ ,

where a,b,c,d are rational numbers and $N=n_1+n_2+n_3$.

However , if $N=n_1+n_2+n_3$ is odd , this equation is degenerate giving no information. However when N is odd, using L'Hospital's rule, we get a different equation of the form a' $T(n_1,n_2,n_3)$+b' $T(n_2,n_3,n_1)$+c' $T(n_3,n_1,n_2)$=d'$\pi^{N-1}$, where a',b',c' are rational numbers and d' is a linear combination of

$$\zeta'(-n_3), \zeta'(-n_3-1), \zeta'(-n_3-2), \ldots\ldots\ldots, \zeta'(-n_3-n_2-n_1)$$

with rational coefficients.

Huard, Williams and Zhang [4] have shown that

$$T(a,b,c) = \sum_{i=1}^{a} \binom{a+b-i-1}{a-i} T(i,0,N-i) + \sum_{i=1}^{b} \binom{a+b-i-1}{b-i} T(i,0,N-i)$$

where a,b≥1 and c≥0 are integers and N=a+b+c. In view of this, to evaluate $T(n_1,n_2,n_3)$ for positive integers $n_1,n_2,n_3$, it is enough to evaluate $T(n_1,0,n_2)$, where $n_1,n_2$ are positive integers. However, when one of the complex numbers $s_1,s_2,s_3$, say $s_3$=0 so that N=$s_1$+$s_2$+0=$s_1$+$s_2$, each of the four equations

$a_iT_1+b_iT_2+c_iT_3=d_i$ for i=1,2,3,4, takes the form

$a(T(s_1,0,s_2)+T(s_2,0,s_1))+bT(s_1,s_2,0)=c$, where a,b,c are some functions of complex variables $s_1,s_2,s_3$. (See Proposition 4 below.) That is,

$$a(\zeta(s_1)\cdot\zeta(s_2) - \zeta(s_1+s_2)) + b\zeta(s_1)\cdot\zeta(s_2) = c.$$

Thus we get no information about the value of $T(s_1,0,s_2)$. However when $s_1=n_1$ and $s_2=n_2$ are natural numbers, $s_3$=0 and N=$n_1$+$n_2$+0 is odd, on using L'Hospital's rule, we get an equation of the type $T(n_1,0,n_2)$- $T(n_2,0,n_1)$=f for some complex

number f. This fact along with the equation

T($n_1$,0,$n_2$)+ T($n_2$,0,$n_1$)=ζ($n_1$)ζ($n_2$)-ζ($n_1$+$n_2$)

enables one to evaluate T($n_1$,0,$n_2$), when N=$n_1$+$n_2$ is odd.

Next, we state a few lemmas.

**Lemma 1 :** I) For Re s<0, $\zeta^{(r)}(s,\alpha)$ is a continuous function of the complex variable α and we have $\zeta^{(r)}(s,0)=\zeta^{(r)}(s)=\zeta^{(r)}(s,1)$ for r≥0.

2) We have for Re s<1, $\int_0^1 \zeta^{(r)}(s,\alpha)d\alpha = 0$

Note : For the proof, see Author [7].

**Lemma 2 :** $\int_0^1 \zeta(s_1,\alpha)\zeta(s_2,\alpha)d\alpha$ is an analytic function

I) of $s_1$ in the region Re $s_1$<0.

II) of $s_2$ in the region Re $s_2$<0.

III) $\int_0^1 \zeta(s_1,\alpha)\zeta(s_2,\alpha)d\alpha$ is a differentiable function of $s_1$ for $s_1$≤0 ; and is a differentiable function of $s_2$ for $s_2$≤0.

Remark : A corresponding result holds in respect of $\int_0^1 \zeta(s_1,\alpha)\zeta(s_2,\alpha)\zeta(s_3,\alpha)d\alpha$.

Proof : The result follows from the fact that $\zeta(s,\alpha)$ is an analytic function of s for each α with $0 \leq \alpha \leq 1$, in the region Re s<0 and $\zeta(s,\alpha)$ is a continuous function of α on the unit interval, when Re s<0. Also note that $\zeta(s,\alpha)$ is a continuous function of α on the unit interval [0,1] for real s≤0 .

**Lemma 3** : We have for complex numbers $s_1,s_2$,

$$\int_0^1 \zeta(s_1,\alpha)\cdot\zeta(s_2,\alpha)d\alpha = 2(2\pi)^{s_1+s_2-2}\Gamma(1-s_1)\Gamma(1-s_2)\cos\frac{\pi}{2}(s_1-s_2)\cdot\zeta(2-s_1-s_2)$$

$$\int_0^1 \zeta(s_1,\alpha)\zeta(s_2,1-\alpha)d\alpha = -2(2\pi)^{s_1+s_2-2}\Gamma(1-s_1)\Gamma(1-s_2)\cos\frac{\pi}{2}(s_1+s_2)\cdot\zeta(2-s_1-s_2)$$

**Note :** These results are known . The proof of this lemma is implicit in the proof of our Theorem .

**Lemma 4** : Let $r_1, r_2 \geq 0$ be integers . For Re $s_1$, Re $s_2<0$ (or for real $s_1, s_2 \leq 0$) , we have

$$\int_0^1 \frac{\partial^{r_1}}{\partial s_1^{r_1}}\zeta(s_1,\alpha)\cdot\frac{\partial^{r_2}}{\partial s_2^{r_2}}\zeta(s_2,\alpha)d\alpha = \frac{\partial^{r_1}}{\partial s_1^{r_1}}\frac{\partial^{r_2}}{\partial s_2^{r_2}}\int_0^1 \zeta(s_1,\alpha)\cdot\zeta(s_2,\alpha)d\alpha$$

$$= \frac{\partial^{r_1}}{\partial s_1^{r_1}}\frac{\partial^{r_2}}{\partial s_2^{r_2}}\left(2(2\pi)^{s_1+s_2-2}\cdot\Gamma(1-s_1)\Gamma(1-s_2)\cos\frac{\pi}{2}(s_1-s_2)\zeta(2-s_1-s_2)\right)$$

**Lemma 5** : We have

I) $\dfrac{\partial}{\partial\alpha}\zeta(s,\alpha) = -s\zeta(s+1,\alpha)$

II) $\dfrac{\partial}{\partial\alpha}\zeta'(s,\alpha) = -\bigl(\zeta(s+1,\alpha) + s\zeta'(s+1,\alpha)\bigr)$

so that for $s \neq 1$.

1) $\zeta(s,\alpha) = \dfrac{\partial}{\partial\alpha}\left(\dfrac{\zeta(s-1,\alpha)}{1-s}\right)$

2) $\zeta'(s,\alpha) = \dfrac{\partial}{\partial\alpha}\left(\dfrac{\zeta(s-1,\alpha)}{(1-s)^2} + \dfrac{\zeta'(s-1,\alpha)}{1-s}\right).$

**Note :** For the proof , see author [7] .

**Lemma 6** : We have for integers $n_1, n_2, n_3 \geq 1$, $\int_0^1 \zeta(1-n_1,\alpha)\cdot\zeta(1-n_2,\alpha)\cdot\zeta(1-n_3,\alpha)d\alpha$

$$= -\frac{1}{n_1 n_2 n_3}\sum_{i=0}^{n_1}\sum_{j=0}^{n_2}\sum_{k=0}^{n_3}\frac{\binom{n_1}{i}\binom{n_2}{j}\binom{n_3}{k}}{i+j+k+1}B_{n_1-i}\cdot B_{n_2-j}\cdot B_{n_3-k}$$

The integral =0 , if $N=n_1+n_2+n_3$ is odd .

**Proof** : $\int_0^1 \zeta(1-n_1,\alpha)\zeta(1-n_2,\alpha)\cdot\zeta(1-n_3,\alpha)d\alpha = -\frac{1}{n_1 n_2 n_3}\int_0^1 B_{n_1}(\alpha)\cdot B_{n_2}(\alpha)\cdot B_{n_3}(\alpha)d\alpha$

$$= -\frac{1}{n_1 n_2 n_3}\int_0^1 \left(\sum_{i=0}^{n_1}\binom{n_1}{i}B_{n_1-i}\alpha^i\right)\left(\sum_{j=0}^{n_2}\binom{n_2}{j}B_{n_2-j}\alpha^j\right)\left(\sum_{k=0}^{n_3}\binom{n_3}{k}B_{n_3-k}\alpha^k\right)$$

$$= -\frac{1}{n_1 n_2 n_3}\sum_i\sum_j\sum_k \binom{n_1}{i}\binom{n_2}{j}\binom{n_3}{k}B_{n_1-i}\cdot B_{n_2-j}\cdot B_{n_3-k}\int_0^1 \alpha^{i+j+k}d\alpha$$

$$= -\frac{1}{n_1 n_2 n_3}\sum_{i=0}^{n_1}\sum_{j=0}^{n_2}\sum_{k=0}^{n_3}\frac{\binom{n_1}{i}\binom{n_2}{j}\binom{n_3}{k}}{i+j+k+1}B_{n_1-i}\cdot B_{n_2-j}\cdot B_{n_3-k}\;.$$

Note that $\int_0^1 B_{n_1}(\alpha)\cdot B_{n_2}(\alpha)\cdot B_{n_3}(\alpha)d\alpha = \int_0^1 B_{n_1}(1-\alpha)\cdot B_{n_2}(1-\alpha)\cdot B_{n_3}(1-\alpha)d\alpha$

$= (-1)^{n_1+n_2+n_3}\int_0^1 B_{n_1}(\alpha)\cdot B_{n_2}(\alpha)\cdot B_{n_3}(\alpha)d\alpha = 0$, if $N = n_1+n_2+n_3$ is odd .

**Lemma 7** : For integers $n_1, n_2, n_3 \geq 1$, we have $\int_0^1 \zeta(1-n_1,\alpha)\cdot\zeta(1-n_2,\alpha)\cdot\zeta'(1-n_3,\alpha)d\alpha$

Is explicitly computable as a linear combination of

$\zeta'(-n_3), \zeta'(-n_3-1), \zeta'(-n_3-2), \ldots\ldots\ldots, \zeta'(-n_3-n_2-n_1)$ with rational coefficients .

**Note :** See author [7] for proof .

**Lemma 8 :** We have $\int_0^1 \zeta^2(0,\alpha) \cdot \zeta'(0,\alpha) d\alpha = -\zeta'(-2)$.

**Proof :** We have $\int_0^1 \zeta^2(0,\alpha) \cdot \zeta'(0,\alpha) d\alpha = \int_0^1 \left(\frac{1}{2}-\alpha\right)^2 \frac{d}{d\alpha}(\zeta(-1,\alpha)+\zeta'(-1,\alpha))$

$= \left[\left(\frac{1}{2}-\alpha\right)^2 (\zeta(-1,\alpha)+\zeta'(-1,\alpha))\right]_{\alpha=0}^1 + \int_0^1 (\zeta(-1,\alpha)+\zeta'(-1,\alpha)) \cdot 2\left(\frac{1}{2}-\alpha\right) d\alpha$

$= 2\int_0^1 \left(\frac{1}{2}-\alpha\right)(\zeta(-1,\alpha)+\zeta'(-1,\alpha)) d\alpha$

$= \int_0^1 (\zeta(-1,\alpha)+\zeta'(-1,\alpha)) d\alpha - 2\int_0^1 \alpha(\zeta(-1,\alpha)+\zeta'(-1,\alpha)) d\alpha = -2\int_0^1 \alpha(\zeta(-1,\alpha)+\zeta'(-1,\alpha)) d\alpha$

Next, $-2\int_0^1 \alpha \zeta'(-1,\alpha) d\alpha = -\int_0^1 \alpha \frac{d}{d\alpha}\left(\frac{\zeta(-2,\alpha)}{2}+\zeta'(-2,\alpha)\right) d\alpha$

$= -\left[\alpha\left(\frac{\zeta(-2,\alpha)}{2}+\zeta'(-2,\alpha)\right)\right]_{\alpha=0}^1 + \int_0^1 \left(\frac{\zeta(-2,\alpha)}{2}+\zeta'(-2,\alpha)\right) d\alpha$

$= -\left(\frac{\zeta(-2)}{2}+\zeta'(-2)\right) + \int_0^1 \left(\frac{\zeta(-2,\alpha)}{2}+\zeta'(-2,\alpha)\right) d\alpha = -\left(\frac{\zeta(-2)}{2}+\zeta'(-2)\right) = -\zeta'(-2)$.

Next, $-2\int_0^1 \alpha \zeta(-1,\alpha) d\alpha = -2\int_0^1 \alpha \frac{d}{d\alpha}\left(\frac{\zeta(-2,\alpha)}{2}\right) d\alpha = -\int_0^1 \alpha \frac{d}{d\alpha}\zeta(-2,\alpha) d\alpha$

$= -\int_0^1 \alpha \frac{d}{d\alpha}\zeta(-2,\alpha) d\alpha = -\left\{[\alpha\zeta(-2,\alpha)]_{\alpha=0}^1 - \int_0^1 \zeta(-2,\alpha) d\alpha\right\} = -\zeta(-2)+0 = -\zeta(-2) = 0$.

Thus $\int_0^1 \zeta^2(0,\alpha) \cdot \zeta'(0,\alpha) d\alpha = -\zeta'(-2)$.

**Lemma 9 :** We have for Re s<1 and for 0<α<1 ,

$$\zeta(s,\alpha) = 2^s \pi^{s-1}\Gamma(1-s)\sum_{n\geq 1}\sin(\frac{\pi s}{2}+2\pi n\alpha)n^{s-1} = \Gamma(1-s)\sum_{|n|\geq 1}e^{2\pi in\alpha}(2\pi in)^{s-1}$$

so that $\zeta^{(r)}(s,\alpha) = \sum_{|n|\geq 1}e^{2\pi in\alpha}\cdot\frac{\partial^r}{\partial s^r}\left(\Gamma(1-s)(2\pi in)^{s-1}\right)$ for $r\geq 0$,

where the value of logarithm stands for its principal value .

**Note : See author [5].** The above series for $\zeta^{(r)}(s,\alpha)$ is its Fourier series on the interval [0,1] as a

function of α .

**Proposition 1) :** I) We have for any complex s ,

$$\int_0^1 \zeta^3(1-s,\alpha)d\alpha = 6(2\pi)^{-3s}\cdot\Gamma^3(s)\cdot\cos\frac{\pi s}{2}\cdot T(s,s,s)$$

so that $T(2n,2n,2n) = \frac{(-1)^n(2\pi)^{6n}}{6((2n-1)!)^3}\cdot\int_0^1 B_{2n}^3(\alpha)d\alpha$

and thus T(2n,2n,2n) is a rational multiple of $\pi^{6n}$ for any integer $n\geq 1$.

II) We have for any integer $n\geq 0$,

$$T(2n+1,2n+1,2n+1) = \frac{2(2\pi)^{6n+2}(-1)^n}{((2n)!)^3}\int_0^1 \zeta^2(-2n,\alpha)\zeta'(-2n,\alpha)d\alpha$$

Thus $T(2n+1,2n+1,2n+1)$ is $\pi^{6n+2}$ multiplied by a linear combination of

$\zeta'(-2n-1),\zeta'(-2n-2),\zeta'(-2n-3),\ldots\ldots\ldots\ldots,\zeta'(-6n-2)$ with rational coefficients .

In particular , $T(1,1,1) = 2(2\pi)^2 \cdot \int_0^1 \zeta^2(0,\alpha)\cdot\zeta'(0,\alpha)d\alpha = -8\pi^2\cdot\zeta'(-2) = 2\zeta(3)$ .

**Corollary** : We have I) $T(2n-1,2n,2n+1) = \frac{1}{2}T(2n,2n,2n)$ for $n\geq 1$.

II) $T(2n,2n+1,2n+2) = \frac{1}{2}T(2n+1,2n+1,2n+1)$ for $n\geq 0$.

In particular $T(0,1,2) = \frac{1}{2}T(1,1,1) = \zeta(3)$.

**Proof** : I) From the statement of Theorem , we have on putting $s_1=s_2=s_3=s$ ,

$$\int_0^1 \zeta^3(1-s,\alpha)d\alpha = 6(2\pi)^{-3s} \cdot \Gamma^3(s) \cdot \cos\frac{\pi s}{2} \cdot T(s,s,s)$$

From this , we get the expression for $T(2n,2n,2n)$ for integer $n \geq 1$, as above .

II) For $s = 2n+1$ (for integer $n \geq 0$), the equation

$$\int_0^1 \zeta^3(1-s,\alpha)d\alpha = 6(2\pi)^{-3s} \cdot \Gamma^3(s) \cdot \cos\frac{\pi s}{2} \cdot T(s,s,s) \text{ gives either side =0 .}$$

Hence using L'Hospital's rule,

$$\lim_{s \to 2n+1} \left( \frac{\int_0^1 \zeta^3(1-s,\alpha)d\alpha}{\cos\frac{\pi s}{2}} \right) = 6 \cdot (2\pi)^{-3(2n+1)}((2n)!)^3 \cdot T(2n+1,2n+1,2n+1).$$

Now $\lim_{s \to 2n+1} \left( \frac{\int_0^1 \zeta^3(1-s,\alpha)d\alpha}{\cos\frac{\pi s}{2}} \right) = \frac{\lim_{s \to 2n+1}\left(\frac{\partial}{\partial s}\int_0^1 \zeta^3(1-s,\alpha)d\alpha\right)}{\lim_{s \to 2n+1}(-\frac{\pi}{2}\sin\frac{\pi s}{2})}$

$$\lim_{s \to 2n+1} \left( \frac{\int_0^1 \frac{\partial}{\partial s}\zeta^3(1-s,\alpha)d\alpha}{-\frac{\pi}{2}(-1)^n} \right) = \frac{-\lim_{s \to 2n+1}\int_0^1 3\zeta^2(1-s,\alpha)\cdot \zeta'(1-s,\alpha)d\alpha}{-\frac{\pi}{2}(-1)^n}$$

$$= \frac{6}{\pi}(-1)^n \cdot \int_0^1 \zeta^2(-2n,\alpha)\cdot \zeta'(-2n,\alpha)d\alpha .$$

Thus $T(2n+1,2n+1,2n+1) = \frac{2(2\pi)^{6n+2}(-1)^n}{((2n)!)^3} \cdot \int_0^1 \zeta^2(-2n,\alpha)\cdot \zeta'(-2n,\alpha)d\alpha.$

In particular $T(1,1,1) = 2(2\pi)^2 \cdot \int_0^1 \zeta^2(0,\alpha) \cdot \zeta'(0,\alpha) d\alpha = -8\pi^2 \cdot \zeta'(-2) = 2\zeta(3)$ on using Lemma 6

**Proposition 2** : Let $n_1, n_2, n_3 \geq 1$ be integers and let $N = n_1 + n_2 + n_3$ be odd . Then

$$\int_0^1 \zeta(1-n_1,\alpha)\zeta(1-n_2,\alpha)\zeta'(1-n_3,\alpha)d\alpha = 2^{-N}\pi^{1-N} \cdot \Gamma(n_1)\Gamma(n_2)\Gamma(n_3)(-1)^{\frac{N+1}{2}} \cdot$$

$$\cdot \left\{(-1)^{n_3} T(n_1,n_2,n_3) + (-1)^{n_1-1} T(n_2,n_3,n_1) + (-1)^{n_2-1} T(n_3,n_1,n_2)\right\}$$

**Proof** : We have for integers $n_1, n_2 \geq 1$ and for any complex s and with N=n$_1$+n$_2$+s ,

$$\int_0^1 \zeta(1-n_1,\alpha)\zeta(1-n_2,\alpha)\zeta(1-s,\alpha)d\alpha = 2(2\pi)^{-(n_1+n_2+s)}\Gamma(n_1)\Gamma(n_2)\Gamma(s) \cdot$$

$$\cdot \left\{\cos\frac{\pi}{2}(N-2s) \cdot T(n_1,n_2,s) + (-1)^{n_1}\cos\frac{\pi}{2}N \cdot T(n_2,s,n_1) + (-1)^{n_2}\cos\frac{\pi N}{2} \cdot T(s,n_1,n_2)\right\}$$

Consider $\dfrac{\int_0^1 \zeta(1-n_1,\alpha) \cdot \zeta(1-n_2,\alpha)\zeta(1-s,\alpha)d\alpha}{\cos\frac{\pi}{2}(N-2s) \cdot T(n_1,n_2,s) + (-1)^{n_1}\cos\frac{\pi}{2}N \cdot T(n_2,s,n_1) + (-1)^{n_2}\cos\frac{\pi N}{2} \cdot T(s,n_1,n_2)}$ .

$= 2(2\pi)^{-(n_1+n_2+s)} \cdot \Gamma(n_1)\Gamma(n_2) \cdot \Gamma(s)$ .

As $s \to n_3$ ( where $n_3 \geq 1$ is an integer) ,

the right hand side (rhs) $\to 2(2\pi)^{-(n_1+n_2+n_3)} \cdot \Gamma(n_1)\Gamma(n_2)\Gamma(n_3)$ .

Consider the the left hand side(lhs) as $s \to n_3$ ,where $n_1 + n_2 + n_3$ is odd .

As $s \to n_3$ , the numerator of the lhs $\to \int_0^1 \zeta(1-n_1,\alpha)\zeta(1-n_2,\alpha) \cdot \zeta(1-n_3,\alpha) d\alpha = 0$.

As $s \to n_3$, the denominator of the lhs

$\to \cos\frac{\pi}{2}(N-2n_3) \cdot T(n_1,n_2,n_3) + (-1)^{n_1}\cos\frac{\pi N}{2} \cdot T(n_2,n_3,n_1) + (-1)^{n_2}\cos\frac{\pi N}{2} \cdot T(n_3,n_1,n_2)$

$= \cos\frac{\pi N}{2}\left\{(-1)^{n_3} T(n_1,n_2,n_3) + (-1)^{n_1} \cdot T(n_2,n_3,n_1) + (-1)^{n_2} T(n_3,n_1,n_2)\right\}$=0 ,

as $N=n_1+n_2+n_3$ is odd. Thus the lhs is of the form $\dfrac{0}{0}$ as $s \to n_3$.

Using L'Hospital's rule, we have the lhs

$$\dfrac{\lim\limits_{s \to n_3} \dfrac{\partial}{\partial s} \int_0^1 \zeta(1-n_1,\alpha)\zeta(1-n_2,\alpha)\zeta(1-s,\alpha)d\alpha}{\lim\limits_{s \to n_3} \dfrac{\partial}{\partial s}\{\cos\dfrac{\pi}{2}(N-2s)\cdot T(n_1,n_2,s)+(-1)^{n_1}\cos\dfrac{\pi N}{2}\cdot T(n_2,s,n_1)+(-1)^{n_2}\cos\dfrac{\pi N}{2}\cdot T(s,n_1,n_2)\}}$$

The numerator of the lhs $= \lim\limits_{s \to n_3} \int_0^1 \zeta(1-n_1,\alpha)\zeta(1-n_2,\alpha)\dfrac{\partial}{\partial s}\zeta(1-s,\alpha)d\alpha$

$= -\int_0^1 \zeta(1-n_1,\alpha)\zeta(1-n_2,\alpha)\zeta'(1-n_3,\alpha)d\alpha$.

The denominator of the lhs

$= \left\{\cos\dfrac{\pi}{2}(N-2s)\cdot\dfrac{\partial}{\partial s}T(n_1,n_2,s)+(-1)^{n_1}\cos\dfrac{\pi N}{2}\dfrac{\partial}{\partial s}T(n_2,s,n_1)+(-1)^{n_2}\cos\dfrac{\pi N}{2}\dfrac{\partial}{\partial s}T(s,n_1,n_2)\right\}_{s=n_3}$

$+ \dfrac{\pi}{2}\left\{\sin\dfrac{\pi}{2}(N-2n_3)\cdot T(n_1,n_2,,n_3)-\sin\dfrac{\pi}{2}(N-2n_1)\cdot T(n_2,n_3,n_1)-\sin\dfrac{\pi}{2}(N-2n_2)\cdot T(n_3,n_1,n_2)\right\}$

$= 0 + \dfrac{\pi}{2}\left\{(-1)^{n_3}\sin\dfrac{\pi N}{2}\cdot T(n_1,n_2,n_3)+(-1)^{n_1-1}\cdot\sin\dfrac{\pi N}{2}\cdot T(n_2,n_3,n_1)+(-1)^{n_2-1}\cdot\sin\dfrac{\pi N}{2}\cdot T(n_3,n_1,n_2)\right\}$,

where $N=n_1+n_2+n_3$, which is odd.

Thus the denominator of the lhs

$= \dfrac{\pi}{2}\cdot\sin\dfrac{\pi N}{2}\cdot\left\{(-1)^{n_3}T(n_1,n_2,n_3)+(-1)^{n_1-1}\cdot T(n_2,n_3,n_1)+(-1)^{n_2-1}\cdot T(n_3,n_1,n_2)\right\}$

$= \dfrac{\pi}{2}\cdot(-1)^{\frac{N-1}{2}}\left\{(-1)^{n_3}T(n_1,n_2,n_3)+(-1)^{n_1-1}T(n_2,n_3,n_1)+(-1)^{n_2-1}\cdot T(n_3,n_1,n_2)\right\}$

As $\dfrac{numerator}{denomiator}$ of lhs=rhs , this gives

$$\int_0^1 \zeta(1-n_1,\alpha)\zeta(1-n_2,\alpha)\cdot \zeta'(1-n_3,\alpha)d\alpha = 2^{-N}\cdot \pi^{1-N}\cdot\Gamma(n_1)\Gamma(n_2)\Gamma(n_3)(-1)^{\frac{N+1}{2}}\cdot$$

$$\{(-1)^{n_3}\cdot T(n_1,n_2,n_3)+(-1)^{n_1-1}\cdot T(n_2,n_3,n_1)+(-1)^{n_2-1}T(n_3,n_1,n_2)\}$$

**Proposition 3 :** Let $n_1,n_2 > 1$ be integers with N=n₁+n₂ odd . Then we have

$$(-1)^{n_1}\cdot T(n_2,0,n_1)+(-1)^{n_2}T(n_1,0,n_2) = \frac{2(2\pi)^{N-1}(-1)^{\frac{N+1}{2}}}{(n_1-1)!(n_2-1)!}\cdot$$

$$\cdot\left\{\int_0^1 (n_1-1)\zeta(2-n_1,\alpha)\zeta(1-n_2,\alpha)\cdot\zeta'(0,\alpha)d\alpha + \int_0^1 (n_2-1)\zeta(2-n_2,\alpha)\zeta(1-n_1,\alpha)\zeta'(0,\alpha)d\alpha\right\}+\zeta(n_1)\zeta(n_2).$$

**Proof :** We have from the statement I) of Theorem ,

for integers $n_1,n_2 > 1$ and for real $s > 0$ , and for $N = n_1+n_2+s$ ,

$$\int_0^1 \zeta(1-n_1,\alpha)\cdot\zeta(1-n_2,\alpha)(s\zeta(1-s,\alpha))d\alpha = 2(2\pi)^{-(n_1+n_2+s)}\cdot\Gamma(n_1)\Gamma(n_2)\Gamma(s+1)\cdot$$

$$\cdot\left\{\cos\frac{\pi}{2}(N-2s)\cdot T(n_1,n_2,s)+(-1)^{n_1}\cos\frac{\pi N}{2}\cdot T(n_2,s,n_1)+(-1)^{n_2}\cos\frac{\pi N}{2}\cdot T(s,n_1,n_2)\right\}$$

Next for $n_1,n_2 > 1$ and for s>0 , $\int_0^1 \zeta(1-n_1,\alpha)\cdot\zeta(1-n_2,\alpha)(s\zeta(1-s,\alpha))d\alpha$

$$= \int_0^1 \zeta(1-n_1,\alpha)\zeta(1-n_2,\alpha)\frac{\partial}{\partial\alpha}\zeta(-s,\alpha))d\alpha$$

$$= [\zeta(1-n_1,\alpha)\zeta(1-n_2,\alpha)\cdot\zeta(-s,\alpha)]_{\alpha=0}^1 - \int_0^1 \zeta(-s,\alpha)\frac{\partial}{\partial\alpha}(\zeta(1-n_1,\alpha)\zeta(1-n_2,\alpha))d\alpha$$

$$= -\int_0^1 \zeta(-s,\alpha)\frac{\partial}{\partial\alpha}(\zeta(1-n_1,\alpha)\zeta(1-n_2,\alpha))d\alpha$$

Note that $\lim\limits_{s \to 0} \int_0^1 \zeta(-s,\alpha) \frac{\partial}{\partial \alpha}(\zeta(1-n_1,\alpha) \cdot \zeta(1-n_2,\alpha))d\alpha$

$= \int_0^1 \zeta(0,\alpha) \frac{\partial}{\partial \alpha}(\zeta(1-n_1,\alpha)\zeta(1-n_2,\alpha))d\alpha = \int_0^1 \left(\frac{1}{2}-\alpha\right) \frac{\partial}{\partial \alpha}(\zeta(1-n_1,\alpha) \cdot \zeta(1-n_2,\alpha))d\alpha$

$= \left(\frac{1}{2}-\alpha\right)\zeta(1-n_1,\alpha)\zeta(1-n_2,\alpha) \Big|_{\alpha=0}^1 + \int_0^1 \zeta(1-n_1,\alpha)\zeta(1-n_2,\alpha)d\alpha$

$= -\zeta(1-n_1)\zeta(1-n_2) + \int_0^1 \zeta(1-n_1,\alpha)\zeta(1-n_2,\alpha)d\alpha$ =0+0 , if $n_1+n_2$ is odd .

Thus $\lim\limits_{s \to 0} \int_0^1 \zeta(-s,\alpha) \frac{\partial}{\partial \alpha}(\zeta(1-n_1,\alpha) \cdot \zeta(1-n_2,\alpha))d\alpha$ =0 .

Also $\lim\limits_{s \to 0} \left\{ \cos\frac{\pi}{2}(N-2s) \cdot T(n_1,n_2,s) + (-1)^{n_1} \cos\frac{\pi}{2} N \cdot T(n_2,s,n_1) + (-1)^{n_2} \cos\frac{\pi}{2} N \cdot T(s,n_1,n_2) \right\}$

$= \cos\frac{\pi}{2}(n_1+n_2) \cdot \{T(n_1,n_2,0) + (-1)^{n_1} T(n_2,0,n_1) + (-1)^{n_2} T(0,n_1,n_2)\} = 0$, if $n_1+n_2$ is odd .

Consider

$$\lim\limits_{s \to 0} \frac{-\int_0^1 \zeta(-s,\alpha) \frac{\partial}{\partial \alpha}(\zeta(1-n_1,\alpha)\zeta(1-n_2,\alpha))d\alpha}{\left\{ \cos\frac{\pi}{2}(s-n_1-n_2) \cdot T(n_1,n_2,s) + (-1)^{n_1} \cos\frac{\pi}{2} N \cdot T(n_2,s,n_1) + (-1)^{n_2} \cos\frac{\pi}{2} N \cdot T(s,n_1,n_2) \right\}}$$

which is , of the type $\frac{0}{0}$ .

By L'Hospital's rule , we have this limit

$$= \frac{-\lim_{s \to 0} \frac{\partial}{\partial s} \int_0^1 \zeta(-s,\alpha) \cdot \frac{\partial}{\partial \alpha}(\zeta(1-n_1,\alpha)\zeta(1-n_2,\alpha))d\alpha}{\lim_{s \to 0}\{[\cos\frac{\pi}{2}(N-2s)\frac{\partial}{\partial s}T(n_1,n_2,s) + (-1)^{n_1}\cos\frac{\pi}{2}N \cdot \frac{\partial}{\partial s}T(n_2,s,n_1)}$$

$$+ (-1)^{n_2}\cos\frac{\pi}{2}N \cdot \frac{\partial}{\partial s}T(s,n_1,n_2)] - \frac{\pi}{2}[\sin\frac{\pi}{2}(s-n_1-n_2) \cdot T(n_1,n_2,s)$$

$$+ (-1)^{n_1}\sin\frac{\pi}{2}(s+n_1+n_2) \cdot T(n_2,s,n_1) + (-1)^{n_2}\sin\frac{\pi}{2}(s+n_1+n_2) \cdot T(s,n_1,n_2)]\}$$

$$= \frac{\lim_{s \to 0} \frac{\partial}{\partial s} \int_0^1 \zeta(-s,\alpha) \frac{\partial}{\partial \alpha}(\zeta(1-n_1,\alpha)\zeta(1-n_2,\alpha))d\alpha}{\frac{\pi}{2}\{\sin\frac{\pi}{2}(s-n_1-n_2) \cdot T(n_1,n_2,s) + (-1)^{n_1}\sin\frac{\pi}{2}(s+n_1+n_2) \cdot T(n_2,s,n_1)}$$

$$+ (-1)^{n_2}\sin\frac{\pi}{2}(s+n_1+n_2) \cdot T(s,n_1,n_2)\}_{s=0}$$

The denominator

$$= \frac{\pi}{2}\{-\sin\frac{\pi}{2}(n_1+n_2) \cdot T(n_1,n_2,0) + (-1)^{n_1}\sin\frac{\pi}{2}(n_1+n_2) \cdot T(n_2,0,n_1) + (-1)^{n_2}\sin\frac{\pi}{2}(n_1+n_2) \cdot T(0,n_1,n_2)\}$$

$$= \frac{\pi}{2}\sin\frac{\pi}{2}(n_1+n_2)\{-T(n_1,n_2,0) + (-1)^{n_1}T(n_2,0,n_1) + (-1)^{n_2}T(0,n_1,n_2)\}$$

$$= \frac{\pi}{2}(-1)^{\frac{n_1+n_2-1}{2}}(-T(n_1,n_2,0) + (-1)^{n_1}T(n_2,0,n_1) + (-1)^{n_2}T(0,n_1,n_2))$$

**The** numerator

$$= \lim_{s \to 0} \frac{\partial}{\partial s} \int_0^1 \zeta(-s,\alpha) \frac{\partial}{\partial \alpha}(\zeta(1-n_1,\alpha)\zeta(1-n_2,\alpha))d\alpha = -\lim_{s \to 0} \int_0^1 \zeta'(-s,\alpha) \frac{\partial}{\partial \alpha}(\zeta(1-n_1,\alpha)\zeta(1-n_2,\alpha))d\alpha$$

$$= -\int_0^1 \zeta'(0,\alpha) \cdot \frac{\partial}{\partial \alpha}(\zeta(1-n_1,\alpha) \cdot \zeta(1-n_2,\alpha))d\alpha$$

$$\frac{numerator}{denominator} = \frac{\frac{2}{\pi}(-1)^{\frac{n_1+n_2+1}{2}} \int_0^1 \zeta'(0,\alpha) \frac{\partial}{\partial \alpha}(\zeta(1-n_1,\alpha)\zeta(1-n_2,\alpha))d\alpha}{-T(n_1,n_2,0) + (-1)^{n_1} T(n_2,0,n_1) + (-1)^{n_2} T(0,n_1,n_2)}$$

This is equal to $\lim_{s \to 0} 2(2\pi)^{-(n_1+n_2+s)} \cdot \Gamma(n_1)\Gamma(n_2)\Gamma(s+1) = 2(2\pi)^{-(n_1+n_2)} \Gamma(n_1) \cdot \Gamma(n_2)$.

Thus we have $(-1)^{\frac{n_1+n_2+1}{2}} \cdot \frac{2}{\pi} \int_0^1 \zeta'(0,\alpha) \frac{\partial}{\partial \alpha}(\zeta(1-n_1,\alpha)\zeta(1-n_2,\alpha))d\alpha$

$$= 2(2\pi)^{-(n_1+n_2)} \Gamma(n_1)\Gamma(n_2)\{-T(n_1,n_2,0) + (-1)^{n_1} T(n_2,0,n_1) + (-1)^{n_2} T(0,n_1,n_2)\}$$

Thus $-T(n_1,n_2,0) + (-1)^{n_1} T(n_2,0,n_1) + (-1)^{n_2} T(0,n_1,n_2)$

$$= \frac{2^{n_1+n_2} \cdot \pi^{n_1+n_2-1} \cdot (-1)^{\frac{n_1+n_2+1}{2}}}{(n_1-1)!(n_2-1)!} \cdot \int_0^1 \zeta'(0,\alpha) \frac{\partial}{\partial \alpha}(\zeta(1-n_1,\alpha)\zeta(1-n_2,\alpha))d\alpha$$

Thus $(-1)^{n_1} T(n_2,0,n_1) + (-1)^{n_2} T(n_1,0,n_2) = \frac{2(2\pi)^{N-1}(-1)^{\frac{N+1}{2}}}{(n_1-1)!(n_2-1)!} \cdot$

$$\cdot \{(n_1-1)\int_0^1 \zeta'(0,\alpha)\zeta(2-n_1,\alpha)\zeta(1-n_2,\alpha)d\alpha + (n_2-1)\int_0^1 \zeta'(0,\alpha)\zeta(1-n_1,\alpha)\zeta(2-n_2,\alpha)d\alpha\} + \zeta(n_1)\zeta(n_2)$$

**Proposition 4 :** Let $s_1, s_2$ be complex numbers with Re $s_1$, Re $s > 1$ and let

$N = s_1 + s_2 + 0 = s_1 + s_2$. Then we have

$$\zeta(1-s_1) \cdot \zeta(1-s_2)$$

$$= 2(2\pi)^{-N} \Gamma(s_1) \Gamma(s_2) \left\{ \cos\frac{\pi}{2}(s_1+s_2) \cdot T(s_1, s_2, 0) + \cos\frac{\pi}{2}(s_1 - s_2) \cdot [T(s_1, 0, s_2) + T(s_2, 0, s_1) + \zeta(s_1 + s_2)] \right\}$$

**Proof**: From the statement I) of Theorem, with $N = s_1 + s_2 + s_3$, we have

$$\int_0^1 \zeta(1-s_1, \alpha) \cdot \zeta(1-s_2, \alpha)(s_3 \cdot \zeta(1-s_3, \alpha)) d\alpha = 2(2\pi)^{-N} \cdot \Gamma(s_1) \cdot \Gamma(s_2) \cdot \Gamma(s_3 + 1).$$

$$\cdot \{ \cos\frac{\pi}{2}(N - 2s_3) \cdot T(s_1, s_2, s_3) + \cos\frac{\pi}{2}(N - 2s_1) \cdot T(s_2, s_3, s_1) + \cos\frac{\pi}{2}(N - 2s_2) \cdot T(s_3, s_1, s_2) \}.$$

Consider for Re $s_1$, Re $s_2 > 1$ and $s_3 > 0$

$$\int_0^1 \zeta(1-s_1, \alpha) \cdot \zeta(1-s_2, \alpha) \cdot (s_3 \zeta(1-s_3, \alpha)) d\alpha = \int_0^1 \zeta(1-s_1, \alpha) \zeta(1-s_2, \alpha) \frac{\partial}{\partial \alpha} \zeta(-s_3, \alpha) d\alpha$$

$$= \zeta(1-s_1, \alpha) \cdot \zeta(1-s_2, \alpha) \zeta(-s_3, \alpha) \big|_{\alpha=0}^1 - \int_0^1 \zeta(-s_3, \alpha) \frac{\partial}{\partial \alpha} (\zeta(1-s_1, \alpha) \cdot \zeta(1-s_2, \alpha)) d\alpha$$

$$= -\int_0^1 \zeta(-s_3, \alpha) \frac{\partial}{\partial \alpha} (\zeta(1-s_1, \alpha) \cdot \zeta(1-s_2, \alpha)) d\alpha$$

*Letting* $s_3 \to 0$ through real values from the right, we have

$$\lim_{s_3 \to 0} \int_0^1 \zeta(1-s_1, \alpha) \cdot \zeta(1-s_2, \alpha)(s_3 \zeta(1-s_3, \alpha)) d\alpha$$

$$= -\lim_{s_3 \to 0} \int_0^1 \zeta(-s_3, \alpha) \frac{\partial}{\partial \alpha} (\zeta(1-s_1, \alpha) \cdot \zeta(1-s_2, \alpha)) d\alpha$$

$$= -\int_0^1 \zeta(0,\alpha) \frac{\partial}{\partial \alpha}(\zeta(1-s_1,\alpha) \cdot \zeta(1-s_2,\alpha))d\alpha = \int_0^1 (\alpha - \frac{1}{2}) \frac{\partial}{\partial \alpha}(\zeta(1-s_1,\alpha) \cdot (1-s_2,\alpha))d\alpha$$

$$= \left[(\alpha - \frac{1}{2})\zeta(1-s_1,\alpha)\zeta(1-s_2,\alpha)\right]_{\alpha=0}^1 - \int_0^1 \zeta(1-s_1,\alpha)\zeta(1-s_2,\alpha)d\alpha$$

$$= \zeta(1-s_1)\zeta(1-s_2) - \int_0^1 \zeta(1-s_1,\alpha) \cdot \zeta(1-s_2,\alpha)d\alpha$$

$$= \zeta(1-s_1) \cdot \zeta(1-s_2) - 2(2\pi)^{-(s_1+s_2)} \cdot \Gamma(s_1)\Gamma(s_2)\cos\frac{\pi}{2}(s_1-s_2) \cdot \zeta(s_1+s_2).$$

$$Hence, \zeta(1-s_1) \cdot \zeta(1-s_2) - 2(2\pi)^{-(s_1+s_2)} \cdot \Gamma(s_1)\Gamma(s_2)\cos\frac{\pi}{2}(s_1-s_2) \cdot \zeta(s_1+s_2)$$

$$= \lim_{s_3 \to 0} 2(2\pi)^{-N}\Gamma(s_1)\Gamma(s_2)\Gamma(s_3+1).$$

$$\cdot \{\cos\frac{\pi}{2}(N-2s_3) \cdot T(s_1,s_2,s_3) + \cos\frac{\pi}{2}(N-2s_1)T(s_2,s_3,s_1) + \cos\frac{\pi}{2}(N-2s_2)T(s_3,s_1,s_2)\}$$

$$= 2(2\pi)^{-(s_1+s_2)} \cdot \Gamma(s_1)\Gamma(s_2) \cdot \{\cos\frac{\pi}{2}(s_1+s_2) \cdot T(s_1,s_2,0) + \cos\frac{\pi}{2}(s_1-s_2)[T(s_1,0,s_2) + T(s_2,0,s_1)]\}$$

This completes the proof of Proposition 4.

Next we give the proof of Theorem.

**Proof of Theorem**: We shall show that for Re $s_1$, Re $s_2$, Re $s_3<1$,

$$\int_0^1 \zeta(s_1,\alpha) \cdot \zeta(s_2,\alpha)\zeta(s_3,\alpha)d\alpha = 2\Gamma(1-s_1)\Gamma(1-s_2)\Gamma(1-s_3)(2\pi)^{s_1+s_2+s_3-3}.$$

$$\cdot \{\sin\frac{\pi}{2}(s_1+s_2-s_3)T(1-s_1,1-s_2,1-s_3) + \sin\frac{\pi}{2}(s_2+s_3-s_1)T(1-s_2,1-s_3,1-s_1)$$

$$+ \sin\frac{\pi}{2}(s_3+s_1-s_2) \cdot T(1-s_3,1-s_1,1-s_2)\}$$

Replacing $s_i$ by $1-s_i$ for i=1,2,3, we shall get

$$\int_0^1 \zeta(1-s_1,\alpha)\zeta(1-s_2,\alpha)\zeta(1-s_3,\alpha)d\alpha = 2(2\pi)^{-\sum_i s_i} \cdot \Gamma(s_1)\Gamma(s_2)\Gamma(s_3) \cdot$$

$$\cdot \left\{ \cos\frac{\pi}{2}(s_1+s_2-s_3)\cdot T(s_1,s_2,s_3) + \cos\frac{\pi}{2}(s_2+s_3-s_1)\cdot T(s_2,s_3,s_1) + \cos\frac{\pi}{2}(s_3+s_1-s_2)\cdot T(s_2,s_3,s_1) \right\}$$

Hence for Re $s_1$, Re $s_2$, Re $s_3 < 0$, consider $\int_0^1 \zeta(s_1,\alpha)\cdot\zeta(s_2,\alpha)\zeta(s_3,\alpha)d\alpha$.

Note that by Lemma 9, $\zeta(s,\alpha) = \Gamma(1-s)\sum_{|n|\geq 1} e^{2\pi i n\alpha}(2\pi i n)^{s-1}$ for Re $s<1$ and for $0<\alpha<1$, where the value of logarithm stands for its principal value.

Thus $\int_0^1 \zeta(s_1,\alpha)\cdot\zeta(s_2,\alpha)\zeta(s_3,\alpha)d\alpha$

$$= \Gamma(1-s_1)\Gamma(1-s_2)\Gamma(1-s_3)\sum_{\substack{n_1,n_2,n_3 \\ |n_1|,|n_2|,|n_3|\geq 1}} (2\pi i n_1)^{s_1-1}(2\pi i n_2)^{s_2-1}(2\pi i n_3)^{s_3-1}\cdot\int_0^1 e^{2\pi i\alpha(n_1+n_2+n_3)}d\alpha$$

$$= \Gamma(1-s_1)\Gamma(1-s_2)\Gamma(1-s_3)\sum_{\substack{|n_1|, |n_2|,|n_3|\geq 1 \\ n_3=-(n_1+n_2)}} (2\pi i n_1)^{s_1-1}(2\pi i n_2)^{s_2-1}(2\pi i n_3)^{s_3-1}$$

$$= \Gamma(1-s_1)\Gamma(1-s_2)\Gamma(1-s_3)(2\pi)^{\sum_i s_i - 3}\sum_{|n_1|,|n_2|\geq 1}(in_1)^{s_1-1}(in_2)^{s_2-1}(-i(n_1+n_2))^{s_3-1}$$

Consider $\sum_{n_1,n_2\geq 1}(in_1)^{s_1-1}(in_2)^{s_2-1}(-i(n_1+n_2))^{s_3-1}$

$$= i^{s_1-1}\cdot i^{s_2-1}\cdot(-i)^{s_3-1}\cdot\sum_{n_1,n_2\geq 1} n_1^{s_1-1}\cdot n_2^{s_2-1}\cdot(n_1+n_2)^{s_3-1}$$

$$= e^{\frac{\pi i}{2}(s_1-1)}\cdot e^{\frac{\pi i}{2}(s_2-1)}\cdot e^{-\frac{\pi i}{2}(s_3-1)}\cdot\sum_{n_1,n_2\geq 1} n_1^{s_1-1}\cdot n_2^{s-1}\cdot(n_1+n_2)^{s_3-1} = e^{\frac{\pi i}{2}(s_1+s_2-s_3-1)}\cdot T(1-s_1,1-s_2,1-s_3).$$

Next consider $\sum_{n_1,n_2\leq -1}(in_1)^{s_1-1}\cdot(in_2)^{s_2-1}\cdot(-i(n_1+n_2))^{s_3-1}$

$$= \sum_{k_1,k_2\geq 1}(-ik_1)^{s_1-1}\cdot(-ik_2)^{s_2-1}(i(k_1+k_2))^{s_3-1} = (-i)^{s_1-1}\cdot(-i)^{s_2-1}\cdot i^{s_3-1}\sum_{k_1,k_2\geq 1} k_1^{s_1-1}\cdot k_2^{s_2-1}\cdot(k_1+k_2)^{s_3-1}$$

$$= e^{-\frac{\pi i}{2}(s_1-1)-\frac{\pi i}{2}(s_2-1)+\frac{\pi i}{2}(s_3-1)}\cdot T(1-s_1,1-s_2,1-s_3) = e^{-\frac{\pi i}{2}(s_1+s_2-s_3-1)}\cdot T(1-s_1,1-s_2,1-s_3).$$

Thus $\sum\sum_{n_1, n_2 \geq 1}(in_1)^{s_1-1}(in_2)^{s_2-1} \cdot (-i(n_1+n_2))^{s_3-1} + \sum\sum_{n_1, n_2 \leq -1}(in_1)^{s_1-1}(in_2)^{s_2-1} \cdot (-i(n_1+n_2))^{s_3-1}$

$= 2\cos\dfrac{\pi}{2}(s_1+s_2-s_3-1)T(1-s_1,1-s_2,1-s_3) = 2\sin\dfrac{\pi}{2}(s_1+s_2-s_3) \cdot T(1-s_1,1-s_2,1-s_3).$

Next, consider $\sum\sum_{n_1 \geq 1, n_2 \leq -1}(in_1)^{s_1-1}(in_2)^{s_2-1} \cdot (-i(n_1+n_2))^{s_3-1} = \sum\sum_{n_1, k \geq 1}(in_1)^{s_1-1}(-ik)^{s_2-1} \cdot (-i(n_1-k))^{s_3-1}$

$= \left(\sum\sum_{\substack{n_1, k \geq 1 \\ n_1 > k}} + \sum\sum_{\substack{n_1, k \geq 1 \\ k > n_1}}\right)(in_1)^{s_1-1}(-ik)^{s_2-1} \cdot (-i(n_1-k))^{s_3-1} = S_1 + S_2$, say.

Now $S_1 = \sum\sum_{n_1 > k}(in_1)^{s_1-1} \cdot (-ik)^{s_2-1} \cdot (-i(n_1-k))^{s_3-1}.$

Writing $n_1 = k+\ell$, we have $S_1 = \sum\sum_{k, \ell \geq 1}(i(k+\ell))^{s-1}(-ik)^{s_2-1} \cdot (-i\ell)^{s_3-1}$

$= (i)^{s_1-1}(-i)^{s_2-1} \cdot (-i)^{s_3-1}\sum\sum_{k, \ell \geq 1}k^{s_2-1}\ell^{s_3-1}(k+\ell)^{s_1-1} = e^{\frac{\pi i}{2}((s_1-1)-(s_2-1)-(s_3-1))} \cdot T(1-s_2,1-s_3,1-s_1)$

$= e^{-\frac{\pi i}{2}(s_2+s_3-s_1-1)} \cdot T(1-s_2,1-s_3,1-s_1).$

Next $S_2 = \sum\sum_{k > n_1}(in_1)^{s_1-1} \cdot (-ik)^{s_2-1} \cdot (-i(n_1-k))^{s_3-1} = \sum\sum_{k > n_1}(in_1)^{s_1-1} \cdot (-ik)^{s_2-1} \cdot (i(k-n_1))^{s_3-1}$

Writing $k = n_1 + \ell$, we have $S_2 = \sum\sum_{n_1, \ell \geq 1}(in_1)^{s_1-1} \cdot (-i(n_1+\ell))^{s_2-1} \cdot (i\ell)^{s_3-1}$

$= e^{\frac{\pi i}{2}(s_3-1+s_1-1-(s_2-1))} \cdot \sum\sum_{n_1, \ell \geq 1}\ell^{s_3-1} \cdot n_1^{s_1-1} \cdot (n_1+\ell)^{s_2-1} = e^{\frac{\pi i}{2}(s_3+s_1-s_2-1)} \cdot T(1-s_3,1-s_1,1-s_2).$

Thus $S_1 + S_2 = e^{-\frac{\pi i}{2}(s_2+s_3-s_1-1)} \cdot T(1-s_2,1-s_3,1-s_1) + e^{\frac{\pi i}{2}(s_3+s_1-s_2-1)} \cdot T(1-s_3,1-s_1,1-s_2).$

Next consider $\sum\sum_{n_1 \leq -1, n_2 \geq 1}(in_1)^{s_1-1}(in_2)^{s_2-1}(-i(n_1+n_2))^{s_3-1} = \sum\sum_{k, n_2 \geq 1}(-ik)^{s_1-1} \cdot (in_2)^{s_2-1}(-i(n_2-k))^{s_3-1}$

$$= \left( \sum_{n_2 >}\sum_k + \sum_{k>}\sum_{n_2} \right) (-ik)^{s_1-1} \cdot (in_2)^{s_2-1} \left(-i(n_2-k)\right)^{s_3-1} = S_3 + S_4 \text{, say}.$$

Imitating the cases of $S_1$ and $S_2$, we find that $S_3 = e^{-\frac{\pi i}{2}(s_3+s_1-s_2-1)} \cdot T(1-s_3,1-s_1,1-s_2)$

and $S_4 = e^{\frac{\pi i}{2}(s_2+s_3-s_1-1)} \cdot T(1-s_2,1-s_3,1-s_1)$.

Thus $S_1+S_2+S_3+S_4$

$$= 2\cos\frac{\pi}{2}(s_2+s_3-s_1-1)\cdot T(1-s_2,1-s_3,1-s_1) + 2\cos\frac{\pi}{2}(s_3+s_1-s_2-1)\cdot T(1-s_3,1-s_1,1-s_2)$$

$$= 2\sin\frac{\pi}{2}(s_2+s_3-s_1)\cdot T(1-s_2,1-s_3,1-s_1) + 2\sin\frac{\pi}{2}(s_3+s_1-s_2)T(1-s_3,1-s_1,1-s_2)$$

Thus we have for Re $s_1$, Re $s_2$, Re $s_3$<1,

$$\int_0^1 \zeta(s_1,\alpha)\cdot\zeta(s_2,\alpha)\cdot\zeta(s_3,\alpha)d\alpha = 2\Gamma(1-s_1)\Gamma(1-s_2)\Gamma(1-s_3)(2\pi)^{(s_1-1)+(s_2-1)+(s_3-1)}.$$

$$\cdot\{\sin\frac{\pi}{2}(s_1+s_2-s_3)T(1-s_1,1-s_2,1-s_3) + \sin\frac{\pi}{2}(s_2+s_3-s_1)\cdot T(1-s_2,1-s_3,1-s_1)$$

$$+\sin\frac{\pi}{2}(s_3+s_1-s_2)\cdot T(1-s_3,1-s_1,1-s_2)\}$$

This completes the proof of I) of Theorem .

II) Next for Re $s_1$, Re $s_2$, Re $s_3$<1,

$$\int_0^1 \zeta(s_1,\alpha)\zeta(s_2,\alpha)\zeta(s_3,1-\alpha)d\alpha = \Gamma(1-s_1)\Gamma(1-s_2)\Gamma(1-s_3)\cdot$$

$$\cdot\int_0^1 d\alpha \sum_{|n_1|\geq 1} e^{2\pi i n_1\alpha}(2\pi i n_1)^{s_1-1} \sum_{|n_2|\geq 1} e^{2\pi i n_2\alpha}\cdot(2\pi i n_2)^{s_2-1} \sum_{|n_3|\geq 1} e^{2\pi i n_3(1-\alpha)})(2\pi i n_3)^{s_3-1}$$

$$= \Gamma(1-s_1)\Gamma(1-s_2)\Gamma(1-s_3)\cdot\int_0^1 d\alpha \sum_{|n_1|,|n_2|,|n_3|\geq 1} e^{2\pi i \alpha(n_1+n_2-n_3)}\cdot(2\pi i n_1)^{s_1-1}\cdot(2\pi i n_2)^{s_2-1}\cdot(2\pi i n_3)^{s_3-1}$$

$$= \Gamma(1-s_1)\Gamma(1-s_2)\Gamma(1-s_3)(2\pi)^{s_1+s_2+s_3-3} \cdot \sum_{\substack{|n_1|,\ |n_2|,\ |n_3|\geq 1 \\ n_3=n_1+n_2}} \sum \sum (in_1)^{s_1-1}(in_2)^{s_2-1}(in_3)^{s_3-1}.$$

We can proceed as in the proof of I) and get the result as stated in II) of Theorem.

References

[1] O. Espinosa , V.H. Moll , The evaluation of Tornheim double sums Part-I , Journal of Number Theory  116(2006) ,pp 200-229 .

[2] J.G.Huard , K.S. Williams and Z.Y. Zhang , On Tornheim's double series , Acta Arith. 75(1996),pp105-117 .

[3] T. Nakamura , A functional relation for the Tornhiem double zeta function , Acta Arith. 125(2006) , pp 257-263.

[4] T. Nakamura , Double Lerch series and their relations, Aequationes Math.75(2008) ,  pp 251-259 .

[5] V.V. Rane , Instant evaluation and demystification of ζ(n),L(n,χ) that Euler, Ramanujan missed-I (ar Xiv:0801.0884) .

[6] V.V. Rane , Instant evaluation and demystification of ζ(n),L(n,χ) that Euler, Ramanujan missed-II (ar Xiv:0807.2626 ) .

[7] V.V. Rane , The α-calculus-cum-α-analysis of    $\partial^r/\partial s^r \zeta(s,\alpha)$,  (ar Xiv:1107.3521) .

[8]H. Tsumura , The functional relations between the Mordell- Tornheim double zeta functions and Riemann zeta function ,   Math. Proc. Cambridge Phil. Soc. (2007),142, pp 395-